\documentclass[a4paper,12pt]{article}
\usepackage[francais,english]{babel}
\usepackage{latexsym,amssymb}
\usepackage{amsmath}
\usepackage{pb-diagram}
\newtheorem{thm}{Theorem}[section]

\newtheorem{definition}[thm]{Definition}
\newtheorem{pro}[thm]{Proposition}

\def\tr{{\rm {tr}}}

\def\div{{\rm {div}}}

\newcommand{\R}{\mathbb{R}}

\begin{document}
\title{\bf{Energy-Momentum Tensor on Foliations}}
\author{{\small {\bf {Georges Habib}}}\\
{\small Institut {\'E}lie Cartan,
Universit{\'e} Henri Poincar{\'e}, Nancy I, B.P. 239}\\
{\small 54506 Vand\oe uvre-L{\`e}s-Nancy Cedex, France}\\
{\small {\tt habib@iecn.u-nancy.fr}}}
\date{}
\maketitle
\begin{abstract}
\noindent
In this paper, we give a new lower bound for the eigenvalues of the Dirac operator on a compact spin manifold. This estimate is motivated by the fact that in its limiting case a skew-symmetric tensor (see Equation \eqref{eq:16}) appears that can be identified  geometrically with the O'Neill tensor of a Riemannian flow, carrying a transversal parallel spinor. The Heisenberg group which is a fibration over the torus is an example of this case. Sasakian manifolds are also considered as particular examples of Riemannian flows. Finally, we characterize the $3$-dimensional case by a solution of the Dirac equation. 
\end{abstract}
\section{Introduction}
The study of the spectrum of the Dirac operator defined on a spin manifold $M,$ has been intensively investigated since it contains subtle information on the geometry of the manifold. 
In \cite{BFGK}, Th. Friedrich proved that on a compact spin manifold $M$ of dimension $n,$ the first eigenvalue $\lambda$ of $D_M$ satisfies 
\begin{equation}
\lambda^2 \geq \frac{n}{4(n-1)}\inf_M \mathrm{Scal}_M,
\label{eq:258}
\end{equation}
where $\mathrm{Scal}_M$ is the scalar curvature of $M,$ supposed positive. The proof is based on the modification of the Levi-Civita connection of the spinor bundle in the direction of the identity and the use of the Schr\"odinger-Lichnerowicz formula \cite{Lic}. The limiting case of \eqref{eq:258} is characterized by the existence of a special section of the spinor bundle called {\it Killing spinor} satisfying an overdetermined differential equation. The manifold is in that case Einstein. \\\\
Observe that Friedrich's estimate contains no information for manifolds with negative or vanishing scalar curvature. 
Hence the estimate established in \cite{h1} for all manifolds (the scalar curvature could be negative) where the author modified the Levi-Civita connection in the direction of a symmetric tensor and leading to a lower bound of the spinorial Laplacian by the norm squared of this tensor. \\\\
In fact, we suppose that on a spin manifold $M,$ there exists a spinor field $\Psi$ such that it satisfies for all $X\in \Gamma(TM),$ 
\begin{equation}
\nabla^M_X\Psi=-E(X)\cdot\Psi,
\label{eq:963}
\end{equation}
where $E$ is a symmetric $2$-tensor defined on $TM.$ Then with the properties of Clifford multiplication, we see that $E$ is equal to the tensor $T^\Psi,$ called the {\it energy-momentum tensor}, defined on the complement set of zeroes of $\Psi$ for all $X, Y \in \Gamma(TM)$ by
\begin{equation}
 T^\Psi(X,Y)=\frac{1}{2}\Re(X\cdot\nabla_Y^M\Psi+Y\cdot\nabla_X^M\Psi,\frac{\Psi}{|\Psi|^2}).  
\label{eq:15}
\end{equation}
Hence he proved that for any eigenspinor $\Psi$ of $D_M$ associated with the first eigenvalue $\lambda$, we have 
\begin{equation}
\lambda^2 \geq \inf_M (\frac{\mathrm{Scal}_M}{4}+|T^\Psi|^2). 
\label{eq:259}
\end{equation}
The important point is that the set of zeroes of $\Psi$ has a Hausdorff dimension equal to $n-2$ (see \cite{Bare}) and hence its measure is zero. The estimate \eqref{eq:259} improves Friedrich's inequality since by the Cauchy-Schwarz inequality, $|T^\Psi|^2 \geqslant \frac{(\tr (T^\Psi))^2}{n}$ (here $\tr$ denotes the trace). The existence of a spinor field satisfying, for all $X \in \Gamma(TM)$ the equation $\nabla^M_X \Psi=-T^\Psi(X)\cdot \Psi,$ characterizes its limiting case. In this case, it is not easy to describe geometrically such manifolds since the lower bound of \eqref{eq:259} depends on the eigenspinor in question.\\\\ 
The study of Equation \eqref{eq:963} in extrinsic spin geometry is the key point for a natural interpretation of this tensor. If the dimension of $M$ is equal to $2,$ Th. Friedrich \cite{T} proved that the existence of a spinor field $\Psi,$ with constant norm satisfying $D_M\Psi=f\Psi,$ where $f$ is a real function on $M,$ is equivalent to the existence of a pair $(\Psi, E)$ satisfying \eqref{eq:963}, where $E$ is a symmetric tensor of trace $f.$ This also implies that $E$ is a Gauss-Codazzi tensor and the manifold $M$ is locally immersed into the Euclidean space $\R^3$ with a mean curvature equal to $f.$ Here we have the following fact \cite{M}: If $M^n$ is a hypersurface of a manifold $N,$ carrying a parallel spinor, then the energy-momentum tensor appears naturally as the second fundamental form $h$ of the hypersurface. Moreover, if the mean curvature $H$ is a positive constant, then we are in the limiting case of the extrinsic estimate established in \cite{HMZ} and we have 
$$\frac{n^2H^2}{4}=\frac{{\rm Scal}_M}{4}+|T^\Psi|^2=\frac{{\rm Scal}_M}{4}+\frac{|h|^2}{4}.$$
In this paper, we study Equation \eqref{eq:963} in a general case. We assume that on a Riemannian spin manifold $(M, g_M),$ there exists a spinor field $\Psi$ satisfying, for all $X\in \Gamma(TM),$ the equation 
\begin{equation}
\nabla^M_X\Psi=-E(X)\cdot\Psi,
\label{eq:987}
\end{equation}
where $E$ is {\it any} endomorphism of $TM$. By using the properties of Clifford multiplication, we find that the symmetric part of $E$ is $T^\Psi$ and the skew-symmetric part of $E$ is the tensor defined, on the complement set of zeroes of $\Psi,$  by 
\begin{equation}
 Q^\Psi(X,Y)=\frac{1}{2}\Re(Y\cdot\nabla_X^M\Psi-X\cdot\nabla_Y^M\Psi,\frac{\Psi}{|\Psi|^2}),  
\label{eq:16}
\end{equation} 
 for all $X, Y \in \Gamma(TM)$ (see Section \ref{sec:6}). Here is the problem to relate these two tensors to the spectrum of the Dirac operator. We prove that if we modify the Levi-Civita connection in the direction of these two tensors, the spinorial Laplacian is bounded from below by the sum of the norm squared of these two tensors. Thus we have:
 \begin{thm} \label{thm:22} Let $(M,g_M)$ be a compact spin manifold, then the first eigenvalue of the Dirac operator satisfies
\begin{equation}
\lambda^2\geq \inf_M(\frac{\mathrm{Scal}_M}{4}+|T^\Psi|^2+|Q^\Psi|^2),
\label{eq:123}
\end{equation}
where $\Psi$ is an eigenspinor of $D_M^2$ associated with $\lambda^2$. The equality case of \eqref{eq:123} is characterized by a solution of \eqref{eq:987} 
\end{thm} 
The Heisenberg group ${\rm Nil}_3$ and the solvable group ${\rm Sol}_3$ are examples of limiting manifolds with negative scalar curvature (the term $Q^\Psi$ is equal to zero, see Examples $1$ and $2$), also the Riemannian product $\mathbb{S}^1\times \mathbb{S}^2$ with positive scalar curvature (the term $T^\Psi$ is equal to zero, see Example $3$).\\\\ 
The study of foliations and in particular the transverse geometry of {\it Riemannian flows} \cite{C}, which are locally given by Riemannian submersions with $1$-dimensional fibres, will allow for a better understand of the tensor $Q^\Psi.$ In fact, the geometry of the normal bundle $Q$ of a Riemannian flow is completely determined by a natural skew-symmetric tensor, called the {\it O'Neill tensor} \cite{O} (see Equation \eqref{eq:Oneill}), since it is related to the Lie bracket of two sections of $Q$. \\\\
After the identification of the spinor bundles of $M$ and $Q,$ we prove that if the normal bundle carries a parallel spinor $\Psi$, the tensor $Q^\Psi$ plays the role of the  O'Neill tensor (see Proposition \ref{pro:128}). Particular examples of Riemannian flows are provided by {\it Sasakian manifolds} \cite{G1}.
We give necessary conditions on such manifolds for admitting transversal parallel spinors (see Proposition \ref{prop:3}) and we prove that it defines a {\it complex K\"ahlerian Killing spinor} \cite{K1} on the cone constructed over the manifold.\\\\
In the last section, we examine closely the case of $3$-dimensional manifolds. We prove that parallel spinors on the normal bundle correspond to solutions of the Dirac equation on $M,$ with constant norm. Hence we obtain the analogue characterization of surfaces established by Th. Friedrich. \\\\
{\bf Acknowledgment.} The author would like to thank the referee for his valuable comments. He is also grateful to Sebasti\'an Montiel and Oussama Hijazi for their support.   
\section{The estimate}\label{sec:6}
\setcounter{equation}{0}
In this section, we prove Theorem \ref{thm:22}. For this, let $(M^n,g_M)$ be a Riemannian spin manifold and let $\nabla^M$ be the Levi-Civita connection associated with $g_M$. We denote by $\Sigma M$ its spinor bundle and we suppose that there exists a spinor field $\Psi$ which satisfies Equation \eqref{eq:987}. As a first consequence of the existence of such a spinor is that its norm is constant. Moreover, by the fact that for all $Z, W \in \Gamma(TM)$, we have $\Re(Z\cdot\Psi,W\cdot\Psi)=g_M(Z,W)|\Psi|^2.$ Then for all $X, Y \in \Gamma(TM),$ we obtain  
\begin{equation*}
\Re(X\cdot \nabla^M_Y\Psi+Y\cdot \nabla^M_X\Psi,\frac{\Psi}{|\Psi|^2})=g_M(X,E(Y))+g_M(Y,E(X)).   
\end{equation*}
Hence we find that the symmetric part of $E$ is equal to $T^\Psi$ defined by \eqref{eq:15}. On the other hand, we have similarly for all $X, Y \in \Gamma(TM),$ that 
\begin{equation*}
\Re(Y\cdot \nabla^M_X\Psi-X\cdot \nabla^M_Y\Psi,\frac{\Psi}{|\Psi|^2})=g_M(Y,E(X))-g_M(X,E(Y)).   
\end{equation*} 
We deduce that the skew-symmetric part of $E,$ is equal to the tensor $Q^\Psi$ defined by \eqref{eq:16}. Here the following question arize: Find an inequality whose limiting case could be characterized by \eqref{eq:987}? For this, we will modify the Levi-Civita connection on $M$ in the direction of the two tensors and we will show that the spinorial Laplacian is bounded from below by the norm of these two tensors. Indeed,\\\\ 
{\bf Proof of Theorem \ref{thm:22}.} For any spinor field $\Psi\in \Gamma(\Sigma M)$ and $X\in \Gamma(TM),$ we consider on $\Gamma(\Sigma M)$ the modified connection $\widetilde{\nabla}_X\Psi=\nabla^M_X\Psi+E^\Psi(X)\cdot\Psi,$ where the tensor $E^\Psi$ is defined for all $X, Y\in \Gamma(TM),$ by 
$$E^\Psi(X,Y)=T^\Psi(X,Y)+Q^\Psi(X,Y)=\Re(Y\cdot\nabla^M_X\Psi,\frac{\Psi}{|\Psi|^2}).$$ 
Then for any local frame $\{e_i\}_{i=1,\cdots,n}$ of $\Gamma(TM),$ we compute 
\begin{eqnarray*}
|\widetilde{\nabla}\Psi|^2&=&|\nabla^M\Psi|^2+|E^\Psi|^2|\Psi|^2-2\sum_{i=1}^n \Re(E^\Psi(e_i)\cdot\nabla^M_{e_i}\Psi,\Psi)\\
&=&|\nabla^M\Psi|^2+|E^\Psi|^2|\Psi|^2-2\sum_{i,j=1}^n g_M(E^\Psi(e_i),e_j)\Re(e_j\cdot\nabla^M_{e_i}\Psi,\Psi)\\
&=&|\nabla^M\Psi|^2-|E^\Psi|^2|\Psi|^2.
\end{eqnarray*}
We then conclude the estimate with the help of the Schr\"odinger-Lichnerowicz formula and the fact that $|E^\Psi|^2=|T^\Psi|^2+|Q^\Psi|^2,$ since the tensor $T^\Psi$ is symmetric and $Q^\Psi$ is skew-symmetric.  
\hfill$\square$\\\\
As we said before, the estimate \eqref{eq:123} improves Friedrich's inequality for an eigenspinor $\Psi$ of $D_M,$ since we have by the Cauchy-Schwarz inequality 
$$|T^\Psi|^2\geq \frac{({\rm tr}(T^\Psi))^2}{n}=\frac{(\Re(D_M\Psi,\frac{\Psi}{|\Psi|^2}))^2}{n}.$$ 
Now, we will prove an analogue of this inequality for the tensor $Q^\Psi.$ For this, we suppose that $M$ carries a K\"ahler structure and let $J$ be its complex structure. It is well-known that on such manifolds there exists a natural operator defined, for all $\Psi \in \Gamma(\Sigma M)$ by, ${\widetilde D}_M\Psi=\sum_{i=1}^n J(e_i)\cdot \nabla^M_{e_i}\Psi$ \cite{k, k2}. This operator is a self-adjoint operator with respect to the $L^2$-product and has a discrete spectrum, if $M$ is compact. Moreover, we can easily prove that $\widetilde D_M^2=D_M^2$ and it anticommutes with $D_M.$  
Now, we write for all $\Psi\in \Gamma(\Sigma M),$ 
\begin{eqnarray*}
\widetilde D_M\Psi&=&\sum_{i=1}^n J(e_i)\cdot \nabla^M_{e_i}\Psi
=\sum_{i,j=1}^n g_M(J(e_i),e_j)e_j \cdot \nabla^M_{e_i}\Psi\\
&=&\sum_{i<j} g_M(J(e_i),e_j)e_j \cdot \nabla^M_{e_i}\Psi+\sum_{i<j} g_M(J(e_j),e_i)e_i \cdot \nabla^M_{e_j}\Psi\\
&=&\sum_{i<j} g_M(J(e_i),e_j)(e_j \cdot \nabla^M_{e_i}\Psi-e_i \cdot \nabla^M_{e_j}\Psi).
\end{eqnarray*}   
By taking the real part of the hermitian product with $\Psi,$ we find 
\begin{equation}
\Re(\widetilde D_M\Psi,\Psi)=2\sum_{i<j} g_M(J(e_i),e_j)Q^\Psi(e_i,e_j)|\Psi|^2=(J,Q^\Psi)|\Psi|^2.
\label{eq:045}
\end{equation}
Hence by using the Cauchy-Schwarz inequality, we deduce that 
$$|Q^\Psi|^2\geq\frac{|(J,Q^\Psi)|^2}{n}=\frac{(\Re(\widetilde D_M\Psi,\frac{\Psi}{|\Psi|^2}))^2}{n}.$$ 
Then for an eigenspinor $\Psi$ of $\widetilde D_M,$ which corresponds with an eigenspinor of $D_M^2$ and not of $D_M,$ the term $|Q^\Psi|^2$ is bounded from below by $\frac{\lambda^2}{n}$ and the inequality \eqref{eq:123} improves Friedrich's estimate. 
\hfill$\square$\\\\ 
{\bf Remark.} It is well-known that on K\"ahler manifolds, a sharp estimate is established by K.D. Kirchberg in \cite{k, k2} depending on the complex dimension. In fact, we establish in \cite{HG} a new estimate on such manifolds involving the two tensors $T^\Psi$ and $Q^\Psi$ which improves Kirchberg's inequalities.
\section{Case of hypersurfaces}\label{sec:7}
\setcounter{equation}{0}
In the following two sections, we will give a geometric interpretation for the tensors $T^\Psi$ and $Q^\Psi.$ We will see that $T^\Psi$ plays the role of the second fundamental form on a manifold foliated by hypersurfaces while the tensor $Q^\Psi$ plays the role of the O'Neill tensor in the case of Riemannian flows. \\\\
Let $(M,g_M,\mathcal{F})$ be a Riemannian spin manifold of dimension $n+1$ and let $\mathcal{F}$ be a foliation of dimension $n,$ i.e. the vector bundle $L$ on $M$ of tangent vectors to the leaves is of rank $n.$ For all $X\in \Gamma(TM)$ and  $Y\in\Gamma(L),$ we set $\nabla^L_X Y=\pi^\perp(\nabla^M_X Y)$ where $\nabla^M$ is the Levi-Civita connection on  $M$ and $\pi^\perp:TM \longrightarrow L$ is the projection. The connection $\nabla^L$ is a metric connection on $L$ with respect to the induced metric on $M$. We assume that the normal bundle is trivial, that means it is generated by a unit vector field $\nu$. Since $TM=L\oplus \R\nu,$ the bundle $L$ is spin as a vector bundle (see \cite{B, M} for details) and carries a spinor bundle denoted by $\Sigma L.$ The two spinor bundles $\Sigma M$ and $\Sigma L$ are identified by a unitary isomorphism for $n$ even whereas the bundle $\Sigma M$ is identified with two copies of $\Sigma L$ for $n$ odd. If we denote by $\Psi^*$ the spinor field of $\Sigma L$ associated with $\Psi$ by the isomorphism, the Clifford multiplication on $L$ is identified with the one on $M$ for all $Y\in\Gamma(L)$ by
\begin{equation}
Y\cdot_L \Psi^*=(\nu \cdot Y \cdot \Psi)^*.
\label{eq:355}
\end{equation}
{\bf A spinorial Gauss formula}
The connection defined in the previous paragraph allows us to establish the spinorial Gauss formula. For this, we set for all $X\in \Gamma(TM), h(X)=-\nabla^M_X\nu.$ The restriction of $h$ to $L$ is the Weingarten map which is symmetric and we have
\begin{equation}
\nabla^M_X Y=\nabla^L_X Y+g_M(h(X),Y)\nu,
\label{eq:14}
\end{equation}
where $Y\in\Gamma(L).$ Then for all $X\in \Gamma(TM)$ and $\Psi\in\Gamma(\Sigma M),$ we have the Gauss formula \cite{B,tr}
\begin{equation*}
\nabla^{M}_{X}\Psi=\nabla_{X}^L\Psi+\frac{1}{2}h(X)\cdot\nu\cdot \Psi.
\end{equation*}
We recall that the energy-momentum tensor is given for all $X, Y\in \Gamma(TM),$ by 
\begin{equation*}
T^\Phi(X,Y)=\frac{1}{2}\Re(\pi^\perp(X) \cdot_L\nabla_{Y}^L\Phi+\pi^\perp(Y) \cdot_L\nabla_{X}^L\Phi,\frac{\Phi}{|\Phi|^2}),
\end{equation*}
where $\Phi$ is a spinor field in $\Gamma(\Sigma L).$ Now we have the following proposition (see also \cite{M}):
\begin{pro}\label{pro:1}
Let $(M,g_M,\mathcal{F})$ be a Riemannian spin manifold and $\mathcal{F}$ a foliation of codimension $1.$ If $M$ carries a parallel spinor  $\Psi,$ then for all $X,Y\in\Gamma(L),$ we have  
$$T^\Phi(X,Y)=-\frac{1}{2}g_M(h(X),Y)=\frac{1}{4}(\mathcal{L}_{\nu}g_M)(X,Y),$$
where $\Phi=\Psi^*$ and $\mathcal{L}_{\nu}$ is the Lie derivative in the direction of $\nu.$ Moreover the foliation $\mathcal{F}$ is Riemannian (i.e. $h(\nu)=0$) if and only if $T^\Phi(\nu,X)=0$ for all $X\in\Gamma(L).$
\end{pro}
{\bf Proof.} From the Gauss formula and the identification in \eqref{eq:355}, we have for all $X\in \Gamma(TM)$ that
$$\nabla_X^L \Phi=\frac{1}{2}h(X)\cdot_L\Phi.$$
On one hand, for $X,Y\in\Gamma(L),$ we have 
\begin{eqnarray*}
T^\Phi(X,Y)&=&\frac{1}{2}\Re(X \cdot_L\nabla_{Y}^L\Phi+Y \cdot_L\nabla_{X}^L\Phi,\frac{\Phi}{|\Phi|^2})\\
&=&\frac{1}{4}\Re(X\cdot_L h(Y)\cdot_L\Phi+Y\cdot_L h(X)\cdot_L \Phi,\frac{\Phi}{|\Phi|^2})\\
&=&-\frac{1}{4}(g_M(X,h(Y))+g_M(Y,h(X)))=-\frac{1}{2}g_M(h(X),Y).
\end{eqnarray*}
On the other hand, we know that 
$$(\mathcal{L}_{\nu}g_M)(X,Y)=g_M(\nabla^M_X\nu,Y)+g_M(\nabla^M_Y\nu,X)=-2g_M(h(X),Y),$$ 
hence the first part of the proposition. For the second part, we compute for all $X\in\Gamma(L)$
$$
T^\Phi(X,\nu)=\frac{1}{2}\Re(X \cdot_L\nabla_{\nu}^L\Phi,\frac{\Phi}{|\Phi|^2})=\frac{1}{4}\Re(X \cdot_L h(\nu)\cdot_L\Phi,\frac{\Phi}{|\Phi|^2})=-\frac{1}{4}g_M(X,h(\nu)).
$$
The foliation is then Riemannian \cite{p} if and only if $T^\Psi(X,\nu)=0$, and the result is proved. 
\hfill$\square$ 
\section{Case of Riemannian flows}\label{sec:5}
\setcounter{equation}{0}
Now, we consider the case of flows, i.e. the leaves are the integral curves of a vector field defined on the manifold. In this case, the bundle $L$ of tangent vectors is trivial, hence the normal bundle $Q$ will play the role of $L.$ Then submersions will be studied instead  of immersions and more precisely  the study of Riemannian submersions in the case of Riemannian flows. \\\\
For this, let $(M^{n+1},g_M,\mathcal{F})$ be a Riemannian manifold with its Levi-Civita connection $\nabla^M$ and let $\xi$ be the unit vector field that defines the flow $\mathcal{F}.$ 
We denote by $Q$ the normal bundle with its induced metric of $M$ and we consider for $X, Y$ sections of $\Gamma(TM)$ and for $Z, W$ sections of $\Gamma(Q)$ in the rest of the paper. We define a metric connection on $Q$ by $\nabla_X^Q Z=\pi(\nabla^M_X Z)$ where $X\in \Gamma(TM), Z\in\Gamma(Q)$ and $\pi:TM\longrightarrow Q$ is the projection. The connection $\nabla^Q$ is related to the connection $\nabla^M,$ for all $X\in \Gamma(TM)$ and $Z\in\Gamma(Q),$ by
$$\nabla^M_X Z=\nabla_X^Q Z-g_M(h(X),Z)\xi,$$
with $h(X)=\nabla^M_X \xi.$ From now on, we assume that $M$ is a spin manifold. The normal bundle is then spin and carries a spin structure induced from the one of $M,$ as the case of the hypersurfaces. The relation between the connections $\nabla^M$ and $\nabla^Q$ could be easily extended on the corresponding spinor bundles and we have
\begin{equation}
\nabla^M_X \Psi=\nabla_X^Q \Psi+\frac{1}{2}\xi\cdot h(X)\cdot\Psi,
\label{eq:13}
\end{equation} 
for all $\Psi\in \Gamma(\Sigma M)$ and $X\in \Gamma(TM).$ For any spinor field $\Psi\in \Gamma(\Sigma M)$ and $X,Y \in \Gamma(TM),$ we denote by $T_\Psi$ (resp. $Q_\Psi$) the symmetric (resp. skew-symmetric) part of the tensor 
$$E_\Psi(X,Y)=\Re(\xi\cdot Y \cdot\nabla^M_X\Psi,\frac{\Psi}{|\Psi|^2}).$$
{\bf Remark.} We should point out that the spectrum of the Dirac operator could be related to the norm of $E_\Psi,$ as in Section \ref{sec:6}. In fact, we can easily prove that 
\begin{equation}
\lambda^2\geq \inf_M (\frac{\mathrm{Scal}_M}{4}+|E_\Psi|^2)\geq \inf_M (\frac{\mathrm{Scal}_M}{4}+|E_\Psi|_Q^2),
\label{eq:estima}
\end{equation}
where $|E_\Psi|_Q^2$ is the norm of $E_\Psi$ evaluated on vectors orthogonal to $\xi.$ We have the following theorem:
\begin{thm} \label{pro:5} Let $(M,g_M,\mathcal{F})$ be a Riemannian spin manifold of dimension $n+1$ and let $\mathcal{F}$ be a flow of $M.$ If the normal bundle admits a parallel spinor $\Phi=\Psi^*,$ then for all $Z, W \in \Gamma(Q),$ we have 
$$T_\Psi(Z,W)=-\frac{1}{4}(\mathcal{L}_\xi g_M)(Z,W) \quad\text{and}\quad Q_\Psi(Z,W)=\frac{1}{4}g_M([Z,W],\xi),$$  
where $\mathcal{L}_\xi$ denotes the Lie derivative in the direction of $\xi.$ Moreover the foliation is minimal (i.e. $\nabla^M_\xi\xi=0$) if and only if $T_\Psi(\xi,Z)=0.$
\end{thm}
{\bf Proof.} If the manifold $M$ admits a transversal parallel spinor $\Phi$, then by \eqref{eq:13} we obtain for all $X\in \Gamma(TM),$ that $\nabla^M_X\Psi =\frac{1}{2}\xi\cdot h(X)\cdot\Psi.$ Hence, for all $Z, W \in \Gamma(Q),$ we deduce  
\begin{eqnarray*}
T_\Psi(Z,W)&=&\frac{1}{2}\Re(\xi\cdot Z \cdot\nabla^M_W\Psi+\xi\cdot W\cdot\nabla^M_Z\Psi,\frac{\Psi}{|\Psi|^2})\\
&=&\frac{1}{4}\Re(\xi\cdot Z \cdot\xi\cdot h(W)\cdot\Psi+\xi\cdot W\cdot\xi\cdot h(Z)\cdot\Psi,\frac{\Psi}{|\Psi|^2})\\
&=&\frac{1}{4}(-g_M(Z,h(W))-g_M(W,h(Z)))=-\frac{1}{4}(\mathcal{L}_\xi g_M)(Z,W).
\end{eqnarray*}
Now we compute 
$$T_\Psi(\xi,Z)=\frac{1}{2}\Re(\xi\cdot\xi\cdot\nabla^M_Z\Psi+\xi\cdot Z\cdot\nabla^M_\xi\Psi,\frac{\Psi}{|\Psi|^2})=-\frac{1}{4}g_M(Z,h(\xi))=-\frac{1}{4}g_M(\kappa,Z)$$ 
where $\kappa=\nabla^M_\xi \xi$ is the mean curvature of $\mathcal{F}.$ Hence the foliation is minimal if and only if $T_\Psi(\xi,Z)=0.$ Similarly, we have 
\begin{eqnarray*}
Q_\Psi(Z,W)&=&\frac{1}{2}\Re(\xi\cdot W \cdot\nabla^M_Z\Psi-\xi\cdot Z\cdot\nabla^M_W\Psi,\frac{\Psi}{|\Psi|^2})\\
&=&\frac{1}{4}\Re(\xi\cdot W \cdot\xi\cdot h(Z)\cdot\Psi-\xi\cdot Z\cdot\xi\cdot h(W)\cdot\Psi,\frac{\Psi}{|\Psi|^2})\\
&=&\frac{1}{4}(-g_M(W,h(Z))+g_M(Z,h(W)))=\frac{1}{4}g_M([Z,W],\xi).
\end{eqnarray*}
The last equality is a consequence of the fact that the torsion on $M$ is zero. 
\hfill$\square$\\\\
Now we consider  a particular case of flows. A flow is called {\it Riemannian} \cite{C} if for all $Z, W \in  \Gamma(Q),$ we have 
\begin{equation}
(\mathcal{L}_\xi g_M)(Z,W)=0.
\label{eq:026}
\end{equation} 
The metric $g_M$ is said {\it bundle-like} in the sence of \cite{R}. This definition is equivalent to the fact that the restriction $h|_Q$ is skew-symmetric. Moreover, there exists on $Q$ a unique metric connection with vanishing torsion \cite{p}, called {\it transversal Levi-Civita connection}, which it is defined for all $X \in \Gamma(TM)$ and $Z\in\Gamma(Q),$ by 
\begin{equation*}
\nabla _{X} Z =
\left\{\begin{array}{ll}
\pi [\xi,Z], &  \textrm {$X=\xi$},\\\\
\pi (\nabla_{X}^{M}Z), & \textrm {$X\perp \xi$}.
\end{array}\right.
\end{equation*}
An important property for the curvature $R^\nabla$ of the normal bundle is that for all $Y,Z\in \Gamma(Q),$ we have $R^\nabla(\xi,Y)Z=0$ \cite{p}. Hence the operator $R^\nabla(Y, Z):\Gamma(Q)\longrightarrow \Gamma(Q)$ is a well-defined endomorphism. The transversal Ricci operator is defined for all $Y\in \Gamma(Q)$ by ${\rm Ric}^\nabla Y=\sum_{i=1}^n R^\nabla(Y,e_i)e_i,$ where $\{e_i\}_{i=1,\cdots,n}$ is a local frame of $\Gamma(Q).$ The transversal scalar curvature ${\rm Scal}^\nabla$ is the trace of the transversal Ricci curvature. Moreover, the connection $\nabla$ is related to $\nabla^M$ through the Gauss-type formula for all $Z, W \in\Gamma(Q),$ by $\nabla^M_Z W=\nabla_Z W-g_M(h(Z),W)\xi$ and 
\begin{eqnarray}
\nabla^M_\xi Z&=&\nabla^M_Z \xi+[\xi,Z]\nonumber\\
&=&h(Z)+\pi([\xi,Z])+g_M([\xi,Z],\xi)\xi\nonumber\\
&=&h(Z)+\nabla_\xi Z+g_M(\nabla^M_\xi Z-\nabla^M_Z \xi,\xi)\xi\nonumber\\
&=&\nabla_\xi Z+h(Z)-\kappa(Z)\xi.
\label{eq:56}
\end{eqnarray}
Also for the scalar curvatures of $Q$ and $M,$ we have \cite{O}
\begin{equation}
\mathrm{Scal}^\nabla=\mathrm{Scal}_M-2\div_Q \kappa+2|\kappa|^2+|h|_Q^2.
\label{eq:22}
\end{equation}
The geometry of the normal bundle is determined by a skew-symmetric tensor, called the O'Neill tensor \cite{O}, defined for all $X, Y\in \Gamma(TM)$ by
\begin{equation}
A_X Y=\pi^\perp(\nabla^M_{\pi(X)}\pi(Y))+\pi(\nabla^M_{\pi(X)}\pi^\perp(Y)).
\label{eq:Oneill}
\end{equation}
Then if $Z\in\Gamma(Q)$ and $Y=\xi,$ we have $A_Z \xi=\pi(\nabla^M_Z \xi)=h(Z).$ Also if $Z, W\in\Gamma(Q),$ then
\begin{equation}
A_Z W=\pi^\perp(\nabla^M_Z W)=g_M(\nabla^M_Z W,\xi)\xi=-g_M(h(Z),W)\xi.
\label{eq:19}
\end{equation} 
Since the map $h|_Q$ is skew-symmetric, the tensor $A$ has also to be skew-symmetric. Then
\begin{equation*}
A_Z W=\pi^\perp(\nabla^M_Z W)=\pi^\perp(\nabla^M_W Z+[Z,W])=A_W Z+\pi^\perp[Z,W],
\end{equation*}
and we deduce that $A_Z W=\frac{1}{2}\pi^\perp[Z,W].$ The bundle $Q$ is then involutive if and only if the tensor $A$ vanishes. In this case, and if the flow is minimal, then by the De Rham decomposition the manifold is locally isometric to a product of manifolds. This product is global if the manifold is complete and simply connected. From now on, we suppose that the manifold $M$ is spin. For all $\Psi\in \Gamma(\Sigma M),$ we have the analogue of the Gauss formula for Riemannian flows,
\begin{equation}
\left\{\begin{array}{ll}
\nabla^M_\xi\Psi=\nabla_\xi\Psi+\frac{1}{4}\sum_{i=1}^n e_i\cdot h(e_i)\cdot\Psi+\frac{1}{2}\xi\cdot\kappa\cdot \Psi, &\textrm {}\\\\
\nabla^M_Z\Psi=\nabla_Z\Psi+\frac{1}{2}\xi\cdot h(Z)\cdot\Psi,&\textrm {}
\end{array}\right.
\label{eq:36}
\end{equation}
where $Z\in\Gamma(Q)$ and $\{e_i\}_{i=1,\cdots,n}$ is a local frame of $\Gamma(Q).$ The proof of the second equality in \eqref{eq:36} is similar to the previous section. For the first one, using Equality \eqref{eq:56}, we write in the frame $\{\xi, e_{1},\cdots,e_{n}\},$
\begin{eqnarray*}
\nabla^{M}_{\xi}\Psi&=& \xi(\Psi)+\frac{1}{2}\sum_{j=1}^n g_M(\nabla^{M}_{\xi}\xi,e_{j})\xi\cdot e_{j}\cdot \Psi+\frac{1}{2}\sum_{i<j}g_M(\nabla^{M}_{\xi}e_{i},e_{j}) e_{i} \cdot e_{j}\cdot \Psi\\ 
&=&\xi(\Psi)+\frac{1}{2}\xi\cdot \nabla^M_\xi \xi\cdot \Psi+\frac{1}{2}\sum_{i<j}g_M(\nabla_\xi e_i+h(e_i),e_{j}) e_{i} \cdot e_{j}\cdot \Psi\\
&=&\nabla_{\xi}\Psi+\frac{1}{2}\xi\cdot \kappa\cdot \Psi+\frac{1}{4}\sum_{i=1}^n e_{i} \cdot h(e_{i})\cdot \Psi. 
\hspace{45mm}\square
\end{eqnarray*}
Now we are ready to state the following proposition:
\begin{pro}\label{pro:128}
Let $(M,g_M,\mathcal{F})$ be a Riemannian spin manifold and let $\mathcal{F}$ be a Riemannian flow. If the normal bundle  carries a parallel spinor $\Phi=\Psi^*,$ then for all $Z, W\in\Gamma(Q),$ we have 
\begin{equation*}
Q_\Psi(Z,W)=\frac{1}{2}g_M(A_Z W,\xi)=-\frac{1}{2}g_M(A_Z\xi,W), 
\end{equation*}
where $A$ denotes the O'Neill tensor.
\end{pro}
{\bf Proof.} From Theorem \ref{pro:5}, we have for all $Z, W\in\Gamma(Q),$
\begin{eqnarray*}
&&Q_\Psi(Z,W)=\frac{1}{4}g_M([Z,W],\xi)=\frac{1}{2}g_M(A_Z W,\xi).\hspace{45mm}\square
\end{eqnarray*}
\section{Case of Sasakian manifolds}
\setcounter{equation}{0}
It is interesting to consider an example of a Riemannian flow. We will discuss the case where the normal bundle admits a parallel spinor. For this, we recall the definition of a Sasakian manifold \cite{G1}. 
\begin{definition}
A Riemannian manifold $(M,g_M)$ of dimension $2m+1$ is called Sasakian, if there exists a unit Killing vector field  $\xi$ such that the tensor $h$ defined for all $X\in \Gamma(TM),$ by $h(X)=\nabla^M_X \xi$ satisfies the following properties: 
\begin{enumerate}
\item $h^2=-{\rm Id}_{TM}+\xi^\flat\otimes \xi,$
\item $(\nabla^M_X h)(Y)=g_M(\xi,Y)X-g_M(X,Y)\xi,$
\end{enumerate}
where $X, Y$ are vector fields in $\Gamma(TM).$ 
\label{def:6}
\end{definition}
Since $\xi$ is a Killing vector field, it then satisfies Equation \eqref{eq:026}. Hence it defines a Riemannian flow with totally geodesic fibres. Moreover, the normal bundle has a K\"ahler structure w.r.t. the connection $\nabla$ defined for all $Z\in \Gamma(Q)$ by $J(Z)=h(Z)$ \cite{G1}. The transversal Ricci tensor is related to the Ricci tensor of $M$ by \cite[eq. 2.5]{G1}
\begin{equation}
\mathrm{Ric}^\nabla Z=\mathrm{Ric}_M Z+2Z \quad \text{and}\quad \mathrm{Ric}_M \xi=2m\xi.
\label{eq:27}
\end{equation}
An important case of Sasakian manifolds is $\eta$-Einstein manifolds (see \cite{G3} and \cite{G2}): 
\begin{definition}
A Sasakian manifold $(M,g_M)$ of dimension $2m+1$ is called $\eta$-Einstein if there exist real functions $\beta$ and $\gamma$ on $M$ such that 
$$\mathrm{Ric}_M=\beta g_M+\gamma\xi^\flat\otimes\xi^\flat.$$
In this case, the functions $\beta $ and $\gamma$ are constant and satisfy $\beta+\gamma=2m.$ The scalar curvature is constant equal to $2m(\beta+1).$
\end{definition}
Let now $(M,g_M)$ be a spin Sasakian manifold. The K\"ahler form $\Omega$ of the bundle $Q$ defined for all $Z,W\in \Gamma(Q)$ by $\Omega(Y,Z):=g_M(J(Y),Z)$ acts on the spinor bundle of $Q$ by \cite{HM4}
$$\Omega=\frac{1}{2}\sum_{i=1}^{2m}e_i\cdot_Q J(e_i)\cdot_Q,$$ 
where $\{e_i\}_{i=1,\cdots,2m}$ is a local frame of $\Gamma(Q).$ It is well-known that under the action of $\Omega$, the spinor bundle of $Q$ splits into an orthogonal sum \cite{HM4}, \cite{k} 
$$\Sigma Q=\oplus_{r=0}^m\Sigma_r Q,$$
where $\Sigma_r Q$ is the eigenbundle of rank $\begin{pmatrix} m\\r \end{pmatrix}$ associated with the eigenvalue $i\mu_r:=i(2r-m)$ of $\Omega.$ Since the bundle $\Sigma M$ is identified with the bundle $\Sigma Q$ by Section \ref{sec:7}, we have the same decomposition for $\Sigma M.$ Moreover, the Killing vector field $\xi$ acts on each eigenbundle $\Sigma_r M$ by \cite{FK}
\begin{equation}
\xi\cdot \Psi_r=(-1)^{r+1}i\Psi_r,
\label{eq:111}
\end{equation}
for all $\Psi_r\in \Gamma(\Sigma_r M).$ Now we have the following proposition:
\begin{pro} \label{prop:3} 
Let $(M,g_M)$ be a simply connected Sasakian spin manifold of dimension $2m+1$ and let $(\xi,h,\eta)$ be its Sasakian structure. If the normal bundle $Q$ admits a parallel spinor $\Phi=\Psi^*,$ then $M$ is $\eta$-Einstein. If moreover the limiting case of Inequality \eqref{eq:estima} is realized, then either $Q$ carries a hyperk\"ahler structure of rank $n=2m=8k$ or the manifold $M$ is isometric to $\R^{4l+1}.$   
\end{pro}
{\bf Proof.} The normal bundle is K\"ahler and Ricci flat with holonomy group is one of the following ${\rm SU}_m, {\rm Sp}_l$ ($m=2l$), $0$ \cite{W}. We deduce by Equation \eqref{eq:22} that $\mathrm{Scal}_M=-2m$ and by \eqref{eq:27} that,
\begin{equation}
\mathrm{Ric}_M Z=-2Z \quad\text{and}\quad \mathrm{Ric}_M\xi = 2m \xi, 
\label{eq:44}
\end{equation}
for all $Z\in\Gamma(Q).$ Thus the manifold $M$ is $\eta$-Einstein and by \eqref{eq:36} it carries a spinor field $\Psi$ that satisfies
\begin{equation}
\left\{\begin{array}{ll}
\nabla^M_\xi\Psi=\frac{1}{2}\Omega\cdot\Psi, &\textrm {}\\\\
\nabla^M_Z\Psi=\frac{1}{2}\xi\cdot h(Z)\cdot\Psi.&\textrm {}
\end{array}\right.
\label{eq:33}
\end{equation}
In this case, the tensor $Q_\Psi(Z,W)=-\frac{1}{2}g_M(h(Z),W)$ for all $Z,W\in\Gamma(Q)$ and we are in the limiting case of Inequality \eqref{eq:estima} if and only if $\Omega\cdot\Psi=0.$ Using the identification of the bundles $\Sigma M$ and $\Sigma Q,$ this condition gives that $\Phi\in\Gamma(\Sigma_{\frac{m}{2}}Q)$ and $m=2l$ is even. Having a holonomy group ${\rm SU}_m,$ the only subbundles that admit parallel spinors in even complex dimension are $\Sigma_0 Q$ and $\Sigma_m Q$ \cite{Mo,Mor,W}, hence a contradiction. Thus the holonomy group is either reduced to ${\rm Sp}_l$ or to $0.$ In the first case the normal bundle admits a hyperk\"ahler structure and the subbundles that admit parallel spinors have the form $\Sigma_{s} Q$ with $s$ even. We deduce that $l$ is even. In the second case, the normal bundle is flat (i.e. $R^\nabla=0$) and by a result of R. Blumenthal \cite[Cor. 2]{Blu}, the manifold $M$ is isometric to $\R^{4l+1}.$
\hfill$\square$\\\\
We recall that a K\"ahler spin manifold $(N^{2m},J,S)$ with complex structure $J$ and spinor bundle $S$ carries a complex K\"ahlerian Killing spinor $\psi=\psi_{r-1}+\psi_{r}\in \Gamma(S_{r-1}\oplus S_r)$ if for each vector field $X$ the differential equations \cite[eq. I.2]{K1}
\begin{equation*}
\left\{\begin{array}{ll}
\nabla^N_X\psi_{r-1}=\frac{c}{2}(X+iJ(X))\cdot\psi_r, &\textrm {}\\\\
\nabla^N_X\psi_r=\frac{c}{2}(X-iJ(X))\cdot\psi_{r-1},&\textrm {}
\end{array}\right.
\end{equation*}
are satisfied, where $c\neq 0$ is a given complex number. Many basic properties have been investigated for a non-trivial solution of the above differential system. In particular, the manifold $N$ is Einstein of odd complex dimension \cite[Thm. 3]{K1}. Now we will relate the particular spinor obtained in Proposition \ref{prop:3} to the cone constucted over the manifold $M$ and we will prove that it corresponds to a complex K\"ahlerian Killing spinor.\\\\ 
For this, let $(M^n,g_M)$ be a Riemannian manifold of dimension $n$ and let $\nabla^M$ be the Levi-Civita connection associated with $g_M.$ We recall the following facts \cite{Ga}, \cite{One}. The cone constructed over $M$ is defined by the Riemannian product $(\mathcal{Z}=\R^+\times M,g_\mathcal{Z}=dt^2\oplus t^2g_M).$ The unit vector field $\partial t$ is orthogonal to the hypersurfaces $M_t=\{t\}\times M\subset\mathcal{Z}$ which foliate the manifold $\mathcal{Z}.$ We denote for all $X\in \Gamma(TM)$ by $h(X)=-\nabla_X^\mathcal{Z}\partial t$ the Weingarten map of $M_t,$ where $\nabla^\mathcal{Z}$ is the Levi-Civita connection associated with $g_\mathcal{Z}.$ We have the following formulas, for all $X,Y\in \Gamma(TM),$  \cite[p. 206]{One}
\begin{center}
$\nabla^\mathcal{Z}_{\partial t}\partial t=0,$\\
$\nabla^\mathcal{Z}_{\partial t}X=\nabla^\mathcal{Z}_{X}\partial t=\frac{1}{t}X,$\\
$\nabla^\mathcal{Z}_{X}Y=\nabla^M_{X}Y-tg_M(X,Y)\partial t.$ 
\end{center}
Using these formulas, we can relate the Ricci curvatures for $M$ and $\mathcal{Z}$ and we have for all $X\in \Gamma(TM),$
$$ {\rm Ric}_\mathcal{Z}\partial t=0,\,\, {\rm Ric}_\mathcal{Z}X=\frac{1}{t^2}({\rm Ric}_M X-(n-1)X),$$
and for the scalar curvatures, we deduce 
$${\rm Scal}_\mathcal{Z}=\frac{1}{t^2}({\rm Scal}_M-n(n-1)).$$
From now on, we suppose that the manifold $M^{2m+1}$ is a Sasakian spin manifold with $(\xi,h,\eta)$ its Sasakian structure. Let the orientation of $\mathcal{Z}$ be such that for any positively orthonormal basis $\{\xi,e_1,\cdots,e_{2m}\}$ of $TM,$ the basis $\{\partial t,\frac{1}{t}\xi,\frac{1}{t}e_1,\cdots,e_{2m}\}$ is positively orthonormal in $\mathcal{Z}.$ Since the dimension of $M$ is odd, the spinor bundles of $M$ and $\mathcal{Z}$ are identified as in Section \ref{sec:7} and we have 
$$\Sigma M\simeq\Sigma \mathcal{Z}^+|_M.$$
Also for the Clifford multiplications, we get from Equation \eqref{eq:355} for all $X\in \Gamma(TM)$ and $\varphi\in \Gamma(\Sigma\mathcal{Z}^+|_M)$ that 
\begin{equation}
X\cdot_M\varphi^*=\frac{1}{t}(\partial t\cdot X\cdot\varphi)^*,
\label{eq:55}
\end{equation} 
where $``\cdot"$ is the Clifford multiplication on $\mathcal{Z}.$ The spinorial Gauss formula is then given for all $X\in \Gamma(TM)$ by 
$$\nabla^\mathcal{Z}_{X}\varphi=\nabla^M_{X}\varphi+\frac{1}{2t}\partial t\cdot X\cdot \varphi,$$
where $\varphi\in \Gamma(\Sigma \mathcal{Z}^+).$
Moreover, we can relate the geometry of $M$ to a particular geometry on the cone. Indeed, the structure $J$ given for all $Y$ orthogonal to $\xi,$ by
$$J(\partial t)=\frac{1}{t}\xi, \,\, J(\xi)=-t\partial t,\,\, J(Y)=h(Y),$$ 
defines a K\"ahler structure on $\mathcal{Z}.$ Let $\Omega^\mathcal{Z}=g_\mathcal{Z}(J(X),Y)$ be the K\"ahler form on the manifold $\mathcal{Z}.$ Its action on the spinor bundle is given by 
\begin{equation}
\Omega^\mathcal{Z}\cdot=\frac{1}{t}\partial t\cdot\xi\cdot+\frac{1}{2t^2}\sum_{i=1}^{2m}e_i\cdot J(e_i)\cdot.
\label{eq:88}
\end{equation}
This formula is a direct consequence from the local expression of $\Omega^\mathcal{Z}$ in the basis $\{\partial t,\frac{1}{t}\xi,\frac{1}{t}e_1,\cdots,\frac{1}{t}e_{2m}\}.$ Now we turn our attention to the cone over the manifolds in Proposition \ref{prop:3}. Using Equations \eqref{eq:44}, we deduce for all $Y\in\Gamma(Q)$ that
$${\rm Ric}_{\mathcal{Z}}\xi={\rm Ric}_{\mathcal{Z}}\partial t=0 \quad\text{and}\quad  {\rm Ric}_{\mathcal{Z}}Y=-\frac{2(m+1)}{t^2}Y.$$
The scalar curvature on $\mathcal{Z}$ is then equal to $-\frac{4m(m+1)}{t^2}.$ Since the spinor field $\Psi=\varphi^*\in \Gamma(\Sigma_{l}Q) (l=\frac{m}{2}$ is supposed even), hence by using Equations \eqref{eq:55} and \eqref{eq:111} we obtain
\begin{equation}
(\frac{1}{t}\partial t\cdot\xi\cdot\varphi)^*=\xi\cdot_M\Psi=-i\Psi=(-i\varphi)^*.
\label{eq:23}
\end{equation} 
Moreover, the action of the last term in \eqref{eq:88} on the spinor field $\varphi$ is zero, since $\Psi$ is the kernel of the K\"ahler form of $\Gamma(Q).$ We then deduce that $\Omega^{\mathcal{Z}}\cdot\varphi=-i\varphi$ and $\varphi\in \Gamma(\Sigma_{l}\mathcal{Z}).$ Therefore, using Equations \eqref{eq:33} and the Gauss formula, we have by Equation \eqref{eq:55} for all $Y\in \Gamma(Q)$ that 
\begin{equation*}
\left\{\begin{array}{ll}
\nabla^\mathcal{Z}_{\partial t}\varphi=0,\\\\
\nabla^\mathcal{Z}_{\xi}\varphi=-\frac{i}{2}\varphi, \\\\
\nabla^\mathcal{Z}_{Y}\varphi=\frac{1}{2t^2}\xi\cdot J(Y)\cdot \varphi+\frac{1}{2t}\partial t\cdot Y\cdot\varphi.
\end{array}\right.
\end{equation*}
The spinor field defined by $\psi:=i\partial t\cdot\varphi$ is a section of the bundle $\Sigma_{l+1}\mathcal{Z}.$ In fact, using Equations \eqref{eq:88} and \eqref{eq:23}, we compute
$$\Omega^\mathcal{Z}\cdot\psi=\frac{i}{t}\partial t\cdot\xi\cdot\partial t\cdot\varphi=\frac{i}{t}\xi\cdot\varphi=\frac{i}{t}(it\partial t\cdot\varphi)=-\partial t\cdot\varphi=i\psi.$$
Hence $\psi\in \Gamma(\Sigma_{l+1}\mathcal{Z})$ and the pair $(\varphi,\psi)\in \Gamma(\Sigma_l \mathcal{Z})\oplus\Gamma(\Sigma_{l+1}\mathcal{Z}).$ Using Equation \eqref{eq:23}, we write for all $Y\in \Gamma(Q),$ 
\begin{eqnarray*}
\nabla^\mathcal{Z}_{Y}\varphi&=&-\frac{1}{2t^2}J(Y)\cdot\xi\cdot \varphi-\frac{1}{2t}Y\cdot\partial t\cdot\varphi\\
&=&-\frac{i}{2t}J(Y)\cdot\partial t\cdot \varphi-\frac{1}{2t}Y\cdot\partial t\cdot\varphi\\
&=&-\frac{1}{2t}J(Y)\cdot\psi+\frac{i}{2t}Y\cdot\psi=\frac{i}{2t}(Y+iJ(Y))\cdot\psi.
\end{eqnarray*}
Similarly, we compute 
\begin{eqnarray*}
\nabla^\mathcal{Z}_{Y}\psi&=&i\nabla^\mathcal{Z}_{Y}\partial t\cdot\varphi+i\partial t\cdot\nabla^\mathcal{Z}_{Y}\varphi\\
&=&\frac{i}{t}Y\cdot\varphi+i\partial t\cdot(-\frac{i}{2t}J(Y)\cdot\partial t\cdot \varphi-\frac{1}{2t}Y\cdot\partial t\cdot\varphi)\\
&=&\frac{i}{2t}Y\cdot\varphi+\frac{1}{2t}J(Y)\cdot\varphi=\frac{i}{2t}(Y-iJ(Y))\cdot\varphi. 
\end{eqnarray*}
The same equations remain true along the vector field $\xi$ with constant $\frac{i}{4t}.$
\hfill$\square$
\section{Case of $3$-dimensional flows}
\setcounter{equation}{0}
Now we will characterize parallel spinors on the normal bundle when the manifold $M$ is of dimension $3.$ We will prove that the existence of such a spinor is equivalent to the existence of a solution of the Dirac equation and we will find the analogy of the characterization for surfaces. For this, we consider a Riemannian spin manifold $(M,g_M,\mathcal{F})$ of dimension $3$ and a Riemannian flow $\mathcal{F},$ supposed minimal, defined by a unit vector field $\xi.$ We recall that the complex volume form 
$$\omega_3=-\xi\cdot e_1\cdot e_2,$$   
acts as the identity on the spinor bundle $\Sigma M,$ where $\{\xi, e_1, e_2\}$ is a local orthonormal frame of $\Gamma(TM).$ Moreover, we have for all $Z\in\Gamma(Q)$
\begin{equation}
Z\cdot_Q \Psi^*=(\xi\cdot Z\cdot\Psi)^* \quad\text{and}\quad (\xi\cdot\Psi)^*=-i{\overline\Psi^*},
\label{eq:28}
\end{equation}
where $\overline\Psi^*=\omega_2\cdot_Q \Psi^*$ and $\omega_2$ is the complex volume form of $\Sigma Q$ defined by $\omega_2=i e_1\cdot_Q e_2.$ Since the map $h(Z)=\nabla^M_Z \xi$ is skew-symmetric, it can be represented by the following matrix 
 $$\begin{pmatrix}
0&-b\\b&0
 \end{pmatrix},$$
where $b:M \longrightarrow \R$ is a function. We have the following theorem: 
\begin{thm}
Let $(M^3,g_M,\mathcal{F})$ be a compact Riemannian manifold and let $\mathcal{F}$ be a minimal Riemannian flow. Then the following properties are equivalent: 
\begin{enumerate}
\item The normal bundle admits a parallel spinor $\Phi=\Psi^*$.
\item The transversal scalar curvature is non-negative and $\Psi$ is a solution of 
\begin{equation}
D_M\Psi=\frac{b}{2}\Psi,
\label{eq:34}
\end{equation}
with $|\Psi|=1$. 
\end{enumerate}
\end{thm}
{\bf Proof.} For $1\Rightarrow 2,$ the first is trivial since the normal bundle is Ricci flat. For the second part, the norm of $\Psi$ is constant by a direct consequence from the equality $X(|\Psi|^2)=2\Re(\nabla_X\Psi,\Psi)$ for all $X\in \Gamma(TM).$ On the other hand, using Equations \eqref{eq:33} and the fact that
\begin{eqnarray*}
\Omega\cdot\Psi&=&\frac{1}{2}(e_1\cdot h(e_1)\cdot \Psi+e_2\cdot h(e_2)\cdot \Psi)\\
&=&\frac{1}{2}(b\, e_1\cdot e_2\cdot \Psi-b\, e_2\cdot e_1\cdot \Psi)=b\, \xi \cdot\Psi.
\end{eqnarray*}
We compute the Dirac operator of $\Psi$ and we find
\begin{eqnarray*}
D_M\Psi&=&-\frac{b}{2}\Psi+\frac{1}{2}e_1\cdot \xi \cdot h(e_1)\cdot\Psi+\frac{1}{2}e_2\cdot \xi \cdot h(e_2)\cdot\Psi\\
&=&-\frac{b}{2}\Psi+\frac{b}{2}e_1\cdot \xi \cdot e_2\cdot\Psi-\frac{b}{2}e_2\cdot \xi \cdot e_1\cdot\Psi=\frac{b}{2}\Psi.
\end{eqnarray*}
For $2\Rightarrow 1,$ we compute first
\begin{eqnarray}
D_M\Psi&=&\xi\cdot\nabla_\xi^M\Psi+e_1\cdot\nabla_{e_1}^M\Psi+e_2\cdot\nabla_{e_2}^M\Psi\nonumber\\
&=&\xi\cdot(\nabla_\xi\Psi+\frac{b}{2}\xi\cdot\Psi)+e_1\cdot(\nabla_{e_1}\Psi+\frac{b}{2}\xi\cdot e_2\cdot\Psi)
\nonumber\\&&+e_2\cdot(\nabla_{e_2}\Psi-\frac{b}{2}\xi\cdot e_1\cdot\Psi)\nonumber\\  
&=&\xi\cdot\nabla_\xi\Psi+e_1\cdot\nabla_{e_1}\Psi+e_2\cdot\nabla_{e_2}\Psi+\frac{b}{2}\Psi.
\label{eq:66}
\end{eqnarray}
Since $\Psi$ satisfies \eqref{eq:34}, we get by \eqref{eq:28} that $D_{tr}\Phi=\nabla_\xi\Phi$ where $\Phi=\Psi^*$ and $D_{tr}$ is the transversal Dirac operator defined for each spinor field $\Phi\in \Gamma(\Sigma Q)$ by \cite[eq. 1.6]{gk2}
$$D_{tr}\Phi=e_1\cdot_Q\nabla_{e_1}\Phi+e_2\cdot_Q\nabla_{e_2}\Phi.$$
Thus we have  
$$\Re(D_{tr}\Phi,\Phi)=\Re(\nabla_\xi\Phi,\Phi)=\frac{1}{2}\xi(|\Phi|^2).$$
The norm of $\Phi$ is being constant, then $\Re(D_{tr}\Phi,\Phi)=0.$ On the other hand, by the fact that for all $Z\in\Gamma(Q),$ we have $R^\nabla(\xi,Z)\Phi=0,$ then 
\begin{eqnarray*}
D_{tr}^2\Phi&=&D_{tr}(\nabla_\xi\Phi)\\
&=&e_1\cdot_Q\nabla_{e_1}\nabla_\xi\Phi+e_2\cdot_Q\nabla_{e_2}\nabla_\xi\Phi\\
&=&e_1\cdot_Q(\nabla_\xi\nabla_{e_1}\Phi+ \nabla_{[e_1,\xi]}\Phi)+e_2\cdot_Q(\nabla_\xi\nabla_{e_2}\Phi+\nabla_{[e_2,\xi]}\Phi).
\end{eqnarray*}
If we choose normal coordinates $\{e_1,e_2\}$ at a point $x$ on $M,$ hence the bracket $[e_i,\xi]_x$ vanishes since the foliation is minimal. Thus, $D_{tr}^2\Phi=\nabla_\xi D_{tr}\Phi$ and 
$$\Re(D_{tr}^2\Phi,\Phi)=\Re(\nabla_\xi D_{tr}\Phi,\Phi)=-(D_{tr}\Phi,\nabla_\xi \Phi)=-|D_{tr}\Phi|^2.$$
The integral over $M,$ gives $D_{tr}\Phi=\nabla_\xi \Phi=0.$ Hence the spinor field $\Phi$ is transversally parallel as a consequence of the Schr\"odinger-Lichnerowicz type formula \cite[eq. 2.1]{gk2} and the fact that the transversal scalar curvature is non-negative. 
\hfill$\square$\\\\
Now we give examples of manifolds in dimension $3$ with negative scalar curvatures, which the limiting case of Inequality \eqref{eq:123} is achieved. We also define a particular Riemannian flow on these manifolds with transversal parallel spinors. \\\\  
{\bf Example 1.} Let $M=\mathrm{Nil}_3$ be the Heisenberg group defined by the quotient of 
$$
G:=\left\{\begin{pmatrix}
1&a&c \\ 0&1&b\\0&0&1
\end{pmatrix}; (a, b, c) \in\R^3 \right\},
$$
by the subgroups $G_k\subset G$ of matrices for which $x,y,z$ are integers divisible by some positive integer $k.$ The Heisenberg group carries a left-invariant metric which has the form \cite{D} 
$$ds^2=dx^2+dy^2+\left(\tau(ydx-xdy)+dz \right)^2,$$ 
where $\tau$ is a non-zero constant real number. We easily verify that the frame $\{e_1, e_2, e_3\}$ defined by
$$e_1=\partial x-\tau y \partial z,\,\,\,  e_2=\partial y+ \tau x\partial z,\,\,\,  e_3=\partial z, $$ 
is an orthonormal frame and satisfies
$$[e_1,e_2]=2\tau e_3, \,\,[e_2,e_3]=0, \,\,[e_1,e_3]=0.$$ 
The Christoffel symbols $ \Gamma_{ij}^k=g(\nabla_{e_i} e_j, e_k)$ are given by
\begin{eqnarray*} 
\Gamma_{12}^3=\Gamma_{23}^1=-\Gamma_{21}^3=\tau,\\
\Gamma_{32}^1=-\Gamma_{31}^2=-\Gamma_{13}^2=\tau.
 \end{eqnarray*} 
The other symbols vanish. The Ricci curvature of $M$ is given by the matrix 
$$\mathrm{Ric}_M=\begin{pmatrix}
-2\tau^2&0&0\\0&-2\tau^2&0\\0&0& 2\tau^2
\end{pmatrix}.$$
The scalar curvature of $M$ is then equal to $-2\tau^2.$ Using the local expression of the covariant derivative of a spinor field \cite{HM4}, the spinor bundle $\Sigma M$ admits a spinor field $\Psi$ which verifies,
$$
\nabla^M_{e_1}\Psi=\frac{1}{2} g_M(\nabla^M_{e_1}e_2, e_3)e_2\cdot e_3\cdot\Psi=\frac{1}{2}\tau e_1\cdot\Psi.
$$
Also, we have that $\nabla^M_{e_2}\Psi=\frac{1}{2}\tau e_2\cdot\Psi$ and $\nabla^M_{e_3}\Psi=-\frac{1}{2}\tau e_3\cdot\Psi.$ Hence the spinor field $\Psi$ is an eigenspinor of the Dirac operator associated with the eigenvalue -$\frac{\tau}{2}.$ Moreover, we compute
$$
T^\Psi(e_1,e_1)=\Re(e_1\cdot\nabla^M_{e_1}\Psi,\frac{\Psi}{|\Psi|^2})=\frac{\tau}{2}\Re(e_1\cdot e_1\cdot\Psi,\frac{\Psi}{|\Psi|^2})=-\frac{\tau}{2}.$$ Similarly, we have that $T^\Psi(e_2,e_2)=-\frac{\tau}{2}$ and $T^\Psi(e_3,e_3)=\frac{\tau}{2}.$ The others are equal to zero and also for $Q^\Psi$. We then deduce that $|T^\Psi|^2=\frac{3\tau^2}{4}$ and we get
$$ \inf_M(\frac{{\rm Scal}_M}{4}+|T^\Psi|^2+|Q^\Psi|^2)=\frac{\tau^2}{4}=\lambda^2.$$
The flow defined by $e_3$ is Riemannian and minimal. In fact, the map $h(Y)=\nabla^M_Y e_3$ is given by
$$h(e_1)=-\tau e_2, \,\,\,h(e_2)=\tau e_1, \,\,\,h(e_3)=0.$$ Using Equations \eqref{eq:36}, we can verify that $\Phi=\Psi^*$ is a transversal parallel spinor. Indeed, we have 
$$\nabla_{e_1}\Phi=(\nabla_{e_1}^M\Psi-\frac{1}{2}e_3\cdot h(e_1)\cdot\Psi)^*=(\frac{1}{2}\tau e_1\cdot\Psi+\frac{1}{2}\tau e_3\cdot e_2\cdot \Psi)^*=0,$$
and $\nabla_{e_2}\Phi=\nabla_{e_3}\Phi=0.$ Hence, we find the result in Proposition \ref{pro:128} by computing 
\begin{eqnarray*} 
Q_\Psi(e_1,e_2)&=&\frac{1}{2}\Re(e_3\cdot e_2\cdot \nabla_{e_1}^M\Psi-e_3\cdot e_1\cdot \nabla_{e_2}^M\Psi,\frac{\Psi}{|\Psi|^2})\\
&=&\frac{\tau}{4}\Re(e_3\cdot e_2\cdot e_1\cdot \Psi -e_3\cdot e_1\cdot e_2\cdot \Psi,\frac{\Psi}{|\Psi|^2})\\
&=&-\frac{\tau}{2}\Re(e_1\cdot e_2 \cdot e_3\cdot\Psi,\frac{\Psi}{|\Psi|^2})=\frac{\tau}{2}=-\frac{1}{2}g_M(h(e_1),e_2).
\hspace{15mm}\square
\end{eqnarray*}
{\bf Example 2.} Let $M$ be the solvable group $\mathrm{Sol}_3.$ The manifold $M$ is the semi-direct product $\R\ltimes \R^2,$ where $t\in \R$ acts on $\R^2$ via the transformation $(x,y)\longrightarrow (e^t x,e^{-t} y).$ We identify $\mathrm {Sol}_3$ with $\R^3$ and the group multiplication is defined by  
$$(x,y,z)\cdot (x',y',z')=(x+e^{-z}x',y+e^z y',z+z').$$
The frame 
$$e_1=e^{-z}\partial x,\,\, e_2=e^z\partial y, \,\, e_3=\partial z,$$
is orthonormal with respect to the left-invariant metric
$$ds^2=e^{2z}dx^2+e^{-2z}dy^2+dz^2.$$ 
We easily verify that the frame $\{e_1,e_2,e_3\}$ satisfies
$$[e_1,e_2]=0,\,\, [e_1,e_3]=e_1,\,\, [e_2,e_3]=-e_2.$$
The Christoffel symbols are given by 
$$ \Gamma_{11}^3=\Gamma_{23}^2=-\Gamma_{13}^1=-\Gamma_{22}^3=-1.$$ 
The other symbols vanish. The scalar curvature is equal to $-2.$ As the previous example, there exists a spinor field $\Psi$ on $\Gamma(\Sigma M)$ which satisfies 
$$\nabla^M_{e_1}\Psi=\frac{1}{2}e_2\cdot \Psi,\,\, \nabla^M_{e_2}\Psi=\frac{1}{2}e_1\cdot \Psi,\,\, \nabla^M_{e_3}\Psi=0.$$
The spinor field $\Psi$ is then a harmonic spinor and we have 
\begin{eqnarray*}
T^\Psi(e_1,e_2)&=&\frac{1}{2}\Re(e_1\cdot\nabla^M_{e_2}\Psi+e_2\cdot\nabla^M_{e_1}\Psi,\frac{\Psi}{|\Psi|^2})\\
&=&\frac{1}{4}\Re(e_1\cdot e_1\cdot\Psi+e_2\cdot e_2\cdot\Psi,\frac{\Psi}{|\Psi|^2})=-\frac{1}{2}.
\end{eqnarray*}
The others are equal to zero and also for $Q^\Psi.$ Hence we deduce that $|T^\Psi|^2=\frac{1}{2}$ and we get 
$$\inf_M(\frac{{\rm Scal}_M}{4}+|T^\Psi|^2+|Q^\Psi|^2)=0=\lambda^2.$$
The flow defined by $e_3$ is minimal and is not Riemannian. In fact, the map $h(Y)=\nabla^M_Y e_3$ satisfies
$$h(e_1)=e_1,\,\, h(e_2)=-e_2,\,\, h(e_3)=0.$$ 
Then, we are in the case of Theorem \ref{pro:5} and we have for all $X\in \Gamma(TM)$ 
$$\nabla^M_X\Psi=\nabla_X\Psi+\frac{1}{2}e_3\cdot h(X)\cdot\Psi.$$
Thus, we find $\nabla_{e_1}\Phi=(\frac{1}{2}e_2\cdot\Psi-\frac{1}{2}e_3\cdot e_1\cdot\Psi)^*=0.$ Also, we deduce that $\nabla_{e_2}\Phi=\nabla_{e_3}\Phi=0$ and $\Phi$ is a parallel spinor on the normal bundle. Now, we compute 
$$T_\Psi(e_1,e_1)=\Re(e_3\cdot e_1\cdot \nabla^M_{e_1}\Psi,\frac{\Psi}{|\Psi|^2})=\frac{1}{2}\Re(e_3\cdot e_1\cdot e_2\cdot \Psi,\frac{\Psi}{|\Psi|^2})=-\frac{1}{2}.$$
On the other hand, we have $-\frac{1}{4}(\mathcal{L}_{e_3}g_M)(e_1,e_1)=-\frac{1}{2}g_M(\nabla^M_{e_1}e_3,e_1)=-\frac{1}{2}.$ Moreover, we write 
\begin{eqnarray*}
Q_\Psi(e_1,e_2)&=&\frac{1}{2}\Re(e_3\cdot e_2\cdot \nabla^M_{e_1}\Psi-e_3\cdot e_1\cdot \nabla^M_{e_2}\Psi,\frac{\Psi}{|\Psi|^2})\\
&=& \frac{1}{4} \Re(e_3\cdot e_2\cdot e_2\cdot\Psi-e_3\cdot e_1\cdot e_1\cdot\Psi,\frac{\Psi}{|\Psi|^2})=0=\frac{1}{4}g_M([e_1,e_2],e_3).
\end{eqnarray*}
\hfill$\square$\\\\
{\bf Example 3.} Let the manifold $M$ be the Riemannian product $\mathbb{S}^1\times \mathbb{S}^2$ and let $\nabla^M$ be the Levi-Civita connection associated with the product metric. The manifold $M$ is a trivial fibration over the sphere $\mathbb{S}^2$ with $\mathbb{S}^1$-fibres. We denote by $\xi$ the unit vector field of the tangent bundle of $\mathbb{S}^1$ and $\{e_1,e_2\}$ is a local orthonormal frame of $\mathbb{S}^2.$ Let $\Phi$ be a Killing spinor on the sphere with Killing number $\frac{1}{2},$ i.e. $\nabla_{e_i}^{\mathbb{S}^2}\Phi=\frac{1}{2}e_i\cdot_{\mathbb{S}^2}\Phi,$ for $i=1,2.$ The scalar curvature on $M$ is then equal to $2.$ Moreover, by using the identification in \eqref{eq:28}, we deduce that the manifold $M$ carries a spinor field $\Psi$ which satisfies 
$$\nabla^M_\xi \Psi=0, \,\,\nabla^M_{e_1} \Psi=\frac{1}{2}e_2\cdot\Psi,\,\,\nabla^M_{e_2} \Psi=-\frac{1}{2}e_1\cdot\Psi.$$  
The spinor field $\Psi$ is an eigenspinor of $D_M^2$ associated with the eigenvalue $1.$ In fact, we have $D_M\Psi=\xi\cdot\Psi$ and 
\begin{eqnarray*}
D_M^2\Psi&=&D_M(\xi\cdot\Psi)=e_1\cdot\nabla^M_{e_1}(\xi\cdot\Psi)+e_2\cdot\nabla^M_{e_2}(\xi\cdot\Psi)+\xi\cdot\nabla^M_\xi(\xi\cdot\Psi)\\
&=&e_1\cdot\xi\cdot\nabla^M_{e_1}\Psi+e_2\cdot\xi\cdot\nabla^M_{e_2}\Psi=-\xi\cdot e_1\cdot e_2\cdot\Psi=\Psi.
\end{eqnarray*}
Moreover, we easily verify that the tensor $T^\Psi$ is equal  to zero and 
\begin{eqnarray*}
Q^\Psi(e_1,e_2)&=&\frac{1}{2}\Re(e_2\cdot\nabla^M_{e_1}\Psi-e_1\cdot\nabla^M_{e_2}\Psi,\frac{\Psi}{|\Psi|^2})\\
&=&\frac{1}{4}\Re(e_2\cdot e_2\cdot\Psi+e_1\cdot e_1\cdot\Psi,\frac{\Psi}{|\Psi|^2})=-\frac{1}{2}. 
\end{eqnarray*}
We also have $Q^\Psi(\xi,e_i)=0$ for $i=1,2.$ Hence we deduce that $|Q^\Psi|^2=\frac{1}{2}$ and 
$$\inf_M(\frac{\mathrm{Scal}_M}{4}+|T^\Psi|^2+|Q^\Psi|^2)=1=\lambda^2.$$
For the Friedrich lower bound, we have $\frac{n}{4(n-1)}\mathop{\mathrm{inf}}\limits_M\mathrm{Scal}_M=\frac{3}{4}.$
\hfill$\square$
\bibliographystyle{amsplain}
\bibliography{Article3}
\end{document}